\documentclass{article}

\newtheorem{theorem}{Theorem}[section]
\newtheorem{definition}{Definition}[section]

\newtheorem{lemma}{Lemma}[section]

\begin{document}
\title{{\bf On Wu and Schaback's Error Bound}}

\author{{\bf Lin-Tian Luh}}

\maketitle


\bigskip


\bigskip

\begin{abstract}
Radial basis functions are a very powerful tool in multivariate approximation since they are mesh free. In the theory of radial basis functions, the most frequently used error bounds are the one raised by Madych and Nelson in \cite{MN2} and the one raised by Wu and Schaback in \cite{WS}, especially the latter. It seems that Wu and Schaback's error bound is five times more frequently used than that of Madych and Nelson's. The reason is that Wu and Schaback's space of approximated functions is much easier to understand, and their error bound is much easier to use. Unfortunately, this error bound contains crucial mistakes which will be pointed out in this paper. Moreover, the r.b.f. people have misunderstood its space of approximands for long. This will also be discussed in depth. These, to some extent, may be a blow to r.b.f. people. However, Madych and Nelson's error bound is still very powerful because it's highly related to Sobolev spaces which contain solutions of many differential equations. All these discussions are based on the central idea of this paper,i.e. unsmooth functions should be approximated by unsmooth functions.
\end{abstract}

\bigskip

\section{Introduction}

\bigskip

Radial basis function interpolation to scattered data $(x_{i}, f_{i}) \in R^{n+1}$ for pairwise distinct points (``centers'') $x_{1}, \ldots , x_{M} \in R^{n}$ uses a function $\phi : R_{\geq 0} \longrightarrow R $ and the space $P_{q}$ of polynomials on $R^{n}$ with order not exceeding $q$ to construct the interpolant

\begin{equation}
  s(x):= \sum _{i=1}^{M} a_{i} \phi (\| x-x_{i}\| ) + \sum _{i=1}^{Q}b_{i}p_{i}(x)    \label{1A}
\end{equation}
via the linear system
\begin{equation}
  \sum_{i=1}^{M}a_{i}\phi (\| x_{j}-x_{i}\| )+ \sum_{i=1}^{Q}b_{i}p_{i}(x_{j})=f_{j}, 1\leq j \leq M    \label{1B}
\end{equation}
$$\sum_{j=1}^{M}a_{j}p_{i}(x_{j})=0,\ \ \ \ \  \ 1\leq i \leq Q=\left( \begin{array}{c}
  q+n-1 \\ n
\end{array}  \right)  $$
where $p_{1},\ldots ,p_{Q}$ are a basis of $P_{q}$, and $f_{j}$ means $f(x_{j})$ with $f$ being the approximated function. For a wide choice of $\phi$ and polynomial orders $q$, including the case $q=Q=0$, the nonsingularity of the $(M+Q)\times (M+Q)$ system (\ref{1B}), written as
\begin{equation}
  \left( \begin{array}{c}
           \begin{array}{cc}
             A & P  \\
             P^{T} & 0 \\
           \end{array}
         \end{array}  \right) \left( \begin{array}{c}
                                       a \\ b 
                                     \end{array} \right) = \left( \begin {array}{c}
                                                                    f \\ 0 
                                                                  \end{array} \right) \label{1C} 
\end{equation}
in matrix notation, has been established by r.b.f. people. A sufficient condition for the linear system to be solvable is that $\psi (x):= \phi (\| x\| ) $ should be conditionally positive definite of order $q$(see \cite{Lu1} for its definition).

\bigskip

The approximated function $f$ should be defined on $\Omega \subseteq R^{n}$. There are no further conditions on $\Omega$. The sampling points $x_{i}$ are allowed to be irregularly distributed over $\Omega$ under the restriction 
\begin{equation}
  p(x_{j})=0,\ 1\leq j \leq M \ for \ p\in P_{q} \ implies \ p=0  \label{1D}
\end{equation}
, which guarantees the uniqueness of the solution of the linear system. For local error estimation, we measure the ``density'' of centers $x_{j}$ near some $x\in \Omega$ by 
$$h_{\rho}(x):=max_{y\in K_{\rho}(x)}min_{1\leq j\leq M}\| y-x_{j}\| $$
for some fixed $\rho > 0$, where $K_{\rho }(x)= \{ y\in R^{n}: \| x-y \| \leq \rho \} $, using the Euclidean norm $\| \cdot \| $.
\section{A brief description of Wu and Schaback's results}
Wu and Schaback showed local error bounds of the following form: Given constants $\rho \in R_{>0},\ q\in N_{\geq 0}$, a radial basis function $\phi$ and a certain function space $F_{\phi}$ to be described below, there exist positive constants $k,m\in N,\ k\geq m,\ h_{0},\ C\in R$ such that for any distribution of points $x_{j}\in R^{n},\ 1\leq j\leq M$, any function $f$ from $F_{\phi}$ and any point $x\in \Omega $ with $h_{\rho}(x)<h_{0}$, the inequality 
\begin{equation}
  |s^{(\mu)}(x)-f^{(\mu)}(x)| \leq c_{f}\cdot C\cdot h_{\rho }^{k-|\mu|}(x)  \label{2A}
\end{equation}
holds for the error and its $\mu$-th derivatives for $0\leq |\mu |\leq m$, where the constant $c_{f}$ depends on $f,\phi$, and $F_{\phi}$ only. Then (\ref{2A}) is called a {\bf local error bound} of order $k$. We define $|\mu|:=\sum _{i}\mu_{i}$ for $\mu\in N^{n}$. 

\bigskip

The linear equations (\ref{1B}) are solvable under the hypotheses of the preceding section. It follows that there is a Lagrange-type representation
\begin{equation}
  s(x)=\sum_{j=1}^{M}f(x_{j})u_{j}(x),\ u_{j}(x_{i})=\delta_{ij},\ 1\leq i,\ j\leq M,  \label{2B}
\end{equation}
of the solution. Introducing vectors
\begin{eqnarray}
  R(x)&:=&(\phi(\| x-x_{1}\|),\ldots , \phi(\| x-x_{M}\| ))^{T} \\
  S(x)&:=&(p_{1}(x),\ldots , p_{Q}(x))^{T},
\end{eqnarray}
we get the following theorem. 
\begin{theorem}
The vector $U(x):=(u_{1}(x),\ldots ,u_{M}(x))^{T}$ formed by the values of the Lagrange basis functions $u_{1},\ldots ,u_{M}$ of (\ref{2B}) at $x\in R^{n}$ coincides with the solution $U_{*}(x)$ of the conditional minimization problem
\end{theorem}
\begin{equation}  
  min\{U^{T}AU-2U^{T}R(x)+\phi(0)\ |\ U\in R^{M},\ P^{T}U=S(x)\}. \label{2C}
\end{equation} \\ 
Proof: omitted.

\bigskip

We can do the same thing for derivatives. If $\phi$ is differentiable of order $|\mu|$ on $(0,\infty )$ and of order $2|\mu|$ around zero, and if $\psi(x):=\phi(\| x\| )$, the solution $U_{*}^{(\mu)}(x)$ of the problem 
\begin{equation}
  min\{ U^{T}AU-2U^{T}R^{(\mu)}(x)+\psi^{(2\mu)}(0)|\ U\in R^{M}, P^{T}U=S^{(\mu)}(x)\} \label{2D}
\end{equation}
uniquely exists. Moreover, the derivative $U^{(\mu)}(x)$ of $U(x)$ coincides with the solution $U_{*}^{(\mu)}(x)$ of (\ref{2D}). It's worth noting that 
\begin{equation}
  \sum_{j=1}^{M}u_{j}^{(\mu)}(x)\cdot p(x_{j})=p^{(\mu)}(x)  \label{2E}
\end{equation}
for all $p\in P_{q}$. The choice of the additive constant in (\ref{2D}) is required for the following definition.
\begin{definition}
  Let $\psi(x)=\phi(\| x\| )$ be conditionally positive definite of order $q$, and assume $\phi\in C^{|\mu|}(0,\infty),\ \phi\in C^{2|\mu|}$ around zero for $\mu\in N_{\geq 0}^{n}$. Then, for any distribution of centers $x_{1},\ldots ,x_{M}$ satisfying (\ref{1D}), the nonnegative function $\kappa _{q}^{(\mu)}$ defined by

 $$ (\kappa _{q}^{(\mu)}(x))^{2}:=min\{ U^{T}AU-2U^{T}R^{(\mu)}(x)+\psi^{(2\mu)}(0)| U\in K_{q}^{(\mu)}(x)\}$$

with the set 
\end{definition}
\begin{equation}
  K_{q}^{(\mu)}:=\{ U=(u_{1},\ldots ,u_{M})^{T}\in R^{M}| \sum_{j=1}^{M}u_{j}p(x_{j})=p^{(\mu)}(x)\ for\ all\ p\in P_{q}\}    \label{2F}
\end{equation}
of admissible vectors is called the {\bf Kriging function} at $x$.

\bigskip

The following theorem is the prototype of Wu and Schaback's main result.
\begin{theorem}
  If the hypotheses of Definition2.1 hold and the Fourier transform $\hat {\psi}$ of $\psi$ satisfies 
\end{theorem}
\begin{equation}
  \psi (y)= \phi (\| y\| )= \frac{1}{(2\pi)^{n}}  \int _{R^{n}}e^{i<y,t>}\hat{\psi}(t)dt,\ y \in R^{n}, \label{2G}
\end{equation}
then
\begin{equation}
  U^{T}AU-2U^{T}R^{(\mu)}(x)+\psi ^{(2\mu)}(0)=\frac{1}{(2\pi )^{n}}\int _{R^{n}} |\sum _{j=1}^{M}u_{j}e^{i<x_{j},t>}-(it)^{\mu}e^{i<x,t>}|^{2}\hat{\psi}(t)dt \label{2H}
\end{equation}
for arbitrary $U=(u_{1},\ldots ,u_{M})\in R^{M}$, expressing the Kriging function via an integral.\\
\\
Proof: The theorem follows from
\begin{eqnarray*}
  &   & \sum_{j,k=1}^{M}w_{j}\overline{w_{k}}\phi (\| x_{j}-x_{k}\| ) \\
  & = & \frac{1}{(2\pi )^{n}}\int _{R^{n}}\sum_{j,k=1}^{M}w_{j}\overline{w_{k}}e^{i<x_{j}-x_{k},t>}\hat{\psi }(t)dt \\
  & = & \frac{1}{(2\pi )^{n}} \int _{R^{n}}|\sum_{j=1}^{M}w_{j}e^{i<x_{j},t>}|^{2}\hat{\psi}(t)dt
\end{eqnarray*}
and
\begin{eqnarray*}
\frac{d^{\mu}}{dx^{\mu}}\phi (\| x-x_{j}\| ) & = & \frac{d^{\mu}}{dx^{\mu}}\frac{1}{(2\pi)^{n}}\int _{R^{n}}e^{i<x-x_{j},t>}\hat{\psi}(t)dt \\
                                             & = & \frac{1}{(2\pi)^{n}}\int _{R^{n}}(it)^{\mu}e^{i<x-x_{j},t>}\hat{\psi}(t)dt 
\end{eqnarray*}

 \hfill $\triangle$ 

However, for most of the interesting radial basis functions $\phi$, (\ref{2G}) does not hold. Hence we must develop an approach of generalized functions.
\begin{theorem}
  Let the generalized Fourier transform of $\psi(x)=\phi(\| x\| )$ exist and coincide with a continuous function $\hat{\psi}$ on $R^{n}-\{0\}$ satisfying
\end{theorem}
\begin{equation}
  0<\hat{\psi}(t)\leq c\left\{ \begin{array}{ll}
                                 \| t\| ^{-n-s_{0}}      & \mbox{for $\| t\| \longrightarrow 0$}\\
                                 \| t\| ^{-n-s_{\infty}} & \mbox{for $\| t\| \longrightarrow \infty $}  
                                \end{array} \right\} \label{2I}
\end{equation}
with constants $c\in R_{>0},\ s_{0},s_{\infty}\in R$, where we additionally assume 
\begin{equation}
  2|\mu| < s_{\infty}\ and \ s_{0} < 2q \label{2J}
\end{equation}
Then for all $U\in K_{q}^{(\mu)}(x)$ we have (\ref{2H}) with a well-defined integral.\\
\\
Proof: Let $g_{U}(t):=\sum_{j=1}^{M}u_{j}e^{i<x_{j},t>}-(it)^{\mu}e^{i<x,t>}$. By the theorem on monotone convergence,
\begin{eqnarray}
  \frac{1}{(2\pi)^{n}} \int _{R^{n}}|g_{U}(t)|^{2}\cdot |\hat{\psi}(t)|dt & = & lim_{m\rightarrow \infty}\frac{1}{(2\pi)^{n}}\int _{R^{n}}|g_{U}(t)|^{2}e^{\frac{-\| t\| ^{2}}{m^{2}}}\hat{\psi}(t)dt  \nonumber \\
                                                                          & = & lim_{m\rightarrow \infty } \frac{1}{(2\pi)^{n}} \int _{R^{n}}\hat{G_{m}}(z)dz 
\end{eqnarray} 
with the test functions 
$$G_{m}(t)=|g_{U}(t)|^{2}e^{\frac{-\| t\| ^{2}}{m^{2}}}$$
using the definition of the generalized Fourier transform for tempered distributions. Furthermore, the Fourier transform $\hat{G}_{m}(z)$ of $G_{m}$ can be explicitly calculated up to a constant $\sigma_{mn}$ as 
\begin{eqnarray*}
  &   & m^{n}\sum_{j,k}u_{j}u_{k}e^{\frac{-\| z-(x_{j}-x_{k})\| ^{2}m^{2}}{4}} \\
  & - & 2m^{n}(-1)^{|\mu|}\sum_{k}u_{k}D^{\mu}(e^{\frac{-\| z-(x-x_{k})\| ^{2}m^{2}}{4}}) \\
  & + & m^{n}D^{2\mu}(e^{\frac{-\| z\| ^{2}m^{2}}{4}}),
\end{eqnarray*}
where $D$ denotes differentiation with respect to $z$. Insertion into (17) yields (\ref{2H}).

\hfill  $\triangle$

We can now write the Kriging function as an integral:
\begin{theorem}
  Under the assumptions of Theorem2.3 the Kriging function has the representation
\end{theorem}
\begin{equation}
  (\kappa_{q}^{(\mu)}(x))^{2}= min_{U\in K_{q}^{(\mu)}(x)}\frac{1}{(2\pi)^{n}}\int _{R^{n}}|\sum_{j=1}^{M}u_{j}e^{i<x_{j},t>}-(it)^{\mu}e^{i<x,t>}|^{2}\hat{\psi}(t)dt \label{2K}
\end{equation}
with $K_{q}^{(\mu)}$ defined as in (\ref{2F}), and the integral exists in the classical sense.\hfill $\triangle$ 
\begin{lemma}
  If the data $(x_{j},f_{j})$ stem from an absolutely integrable real-valued function $f$ on $R^{n}$ with a nicely behaving Fourier transform $\hat{f}$ satisfying
\end{lemma}
\begin{equation}
  f(x)=\frac{1}{(2\pi)^{n}}\int _{R^{n}}e^{i<x,t>}\hat{f}(t)dt,\ x\in R^{n},  \label{2L}
\end{equation}
then
\begin{equation}
  |s^{(\mu)}(x)-f^{(\mu)}(x)|^{2}=|\frac{1}{(2\pi)^{n}}\int _{R^{n}}(\sum_{j=1}^{M}u_{j}^{(\mu)}(x)e^{i<x_{j},t>}-(it)^{\mu}e^{i<x,t>})\hat{f}(t)dt|^{2}  \label{2M}
\end{equation}
where $\phi,\mu$ are as in Definition2.1 and $u_{1}(x),\ldots ,u_{M}(x)$ are Lagrange interpolation functions.\\
\\
Proof: Note that
\begin{eqnarray*}
  s^{(\mu)}(x)-f^{(\mu)}(x) & = & \sum_{j=1}^{M}f(x_{j})u_{j}^{(\mu)}(x)-f^{(\mu)}(x) \\
                            & = & \sum_{j=1}^{M}u_{j}^{(\mu)}(x)\frac{1}{(2\pi)^{n}}\int_{R^{n}}e^{i<x_{j},t>}\cdot \hat{f}(t)dt- \frac{1}{(2\pi)^{n}}\int_{R^{n}}(it)^{\mu}e^{i<x,t>}\hat{f}(t)dt \\
                            & = & \frac{1}{(2\pi)^{n}}\int_{R^{n}}(\sum_{j=1}^{M}u_{j}^{(\mu)}(x)e^{i<x_{j},t>}-(it)^{\mu}e^{i<x,t>})\hat{f}(t)dt.
\end{eqnarray*}
The theorem follows immediately. \hfill  $\triangle$

\bigskip

Assume that 
\begin{equation}
  c_{f}^{2}:=\frac{1}{(2\pi)^{n}}\int _{R^{n}}|\hat{f}(t)|^{2}(\hat{\psi}(t))^{-1}dt < \infty \label{2N}
\end{equation}
in addition to the hypotheses of Definition2.1 and Theorem2.3. Then using the Cauchy-Schwarz inequality, we get 
\begin{eqnarray}
 &   &  |s^{(\mu)}(x)-f^{(\mu)}(x)|^{2} \nonumber  \\
 & = &   \frac{1}{(2\pi)^{n}}\int_{R^{n}}|\sum_{j=1}^{M}u_{j}^{(\mu)}(x)e^{i<x_{j},t>}-(it)^{\mu}e^{i<x,t>}|^{2}\hat{\psi}(t)dt\cdot \frac{1}{(2\pi)^{n}}\int_{R^{n}}|\hat{f}(t)|^{2}(\hat{\psi}(t))^{-1}dt \nonumber  \\
 & = &   (\kappa_{q}^{(\mu)}(x))^{2}\cdot c_{f}^{2}. 
\end{eqnarray} 

\bigskip

To cope with generalized Fourier transforms, we proceed as before and need the following definition.

\begin{definition}
  A function $f:R^{n}\longrightarrow R$ is dominated by a radial basis function $\phi$ satisfying (\ref{2I}) and (\ref{2J}) on $R^{n}-\{ 0\} $, iff $f$ has a generalized Fourier transform $\hat{f}$ coinciding on $R^{n}-\{ 0\} $ with a continuous function satisfying (\ref{2N}) for $\psi(x)=\phi(\| x\| )$.
\end{definition}
Remark: The set $F_{\phi}$ of functions dominated by $\phi$ may be completed to form a Hilbert space with inner product 
$$(f_{1},f_{2})_{\phi}=\int_{R^{n}}\hat{f}_{1}(t)\overline{\hat{f}_{2}(t)}(\hat{\psi}(t))^{-1}dt$$
which was thoroughly studied by Madych and Nelson in \cite{MN1}, \cite{MN2}.
\begin{theorem}
  If $f$ is in the space $F_{\phi}$ of functions dominated by a radial basis function $\phi$ satisfying (\ref{2I})and (\ref{2J}), then the interpolation error can be bounded by
\end{theorem}
\begin{equation}
  |s^{(\mu)}(x)-f^{(\mu)}(x)| \leq \kappa_{q}^{(\mu)}(x)\cdot c_{f} \label{2P}
\end{equation}
, where $c_{f}$ is given by (\ref{2N}).\\
\\
Proof: This is an immediate result of the above arguments. \hfill  $\triangle $

\begin{theorem}
  Let $\phi$ satisfy the assumptions of Theorem2.3, and let $\rho \in R_{>0}$ be given. Then there exist positive real constants $h_{0}$ and $C$ such that for any distribution of centers $x_{i}\in R^{n},\ 1\leq i \leq M$, and any point $x\in R^{n}$ with $h_{\rho}(x)\leq h_{0}$ the Kriging function can be bounded by 
\end{theorem}
\begin{equation}
  \kappa_{q}^{(\mu)}(x)\leq C\cdot h_{\rho}^{\frac{s_{\infty}}{2}-|\mu|}(x)  \label{2Q}
\end{equation}
for $0\leq |\mu| \leq \frac{s_{\infty}}{2}$.\\
\\
Proof: omitted. \hfill  $\triangle $

\bigskip

\section{Crucial mistakes}

The first mistake occurs in (\ref{2H}). The formula holds only when $|\mu| $ is even. Let's check the equation in detail.
\begin{eqnarray*}
  &   & U^{T}AU-2U^{T}R^{(\mu)}(x)+\psi^{(2\mu)}(0) \\
  & = & \frac{1}{(2\pi)^{n}}\int _{R^{n}}|\sum_{j=1}^{M}u_{j}e^{i<x_{j},t>}|^{2}\hat{\psi}(t)dt-2\sum_{j=1}^{M}u_{j}\cdot \frac{1}{(2\pi)^{n}}\int _{R^{n}}(it)^{\mu}e^{i<x-x_{j},t>}\hat{\psi}(t)dt\\
  &   & + \psi^{(2\mu)}(0)\\
  & = & \frac{1}{(2\pi)^{n}}\int _{R^{n}}|\sum_{j=1}^{M}u_{j}e^{i<x_{j},t>}|^{2}\hat{\psi}(t)dt-\frac{1}{(2\pi)^{n}}\int _{R^{n}} 2\sum_{j=1}^{M}u_{j}(it)^{\mu}e^{i<x-x_{j},t>}\hat{\psi}(t)dt\\
  &   & +\psi^{(2\mu)}(0).
\end{eqnarray*}
Now,
\begin{eqnarray*}
  &   & \frac{1}{(2\pi)^{n}}\int _{R^{n}}|\sum_{j=1}^{M}u_{j}e^{i<x_{j},t>}-(it)^{\mu}e^{i<x,t>}|^{2}\hat{\psi}(t)dt\\
  & = & \frac{1}{(2\pi)^{n}}\int _{R^{n}}[\sum_{j=1}^{M}u_{j}e^{i<x_{j},t>}-(it)^{\mu}e^{i<x,t>}][\sum_{j=1}^{M}\overline{u_{j}}e^{-i<x_{j},t>}-\overline{(it)^{\mu}}e^{-i<x,t>}]\\
  &   & \hat{\psi}(t)dt \\
  & = & \frac{1}{(2\pi)^{n}}\int _{R^{n}}\{ |\sum_{j=1}^{M}u_{j}e^{i<x_{j},t>}|^{2}-\sum_{j=1}^{M}u_{j}\overline{(it)^{\mu}}e^{i<x_{j}-x,t>}-\sum_{j=1}^{M}\overline{u_{j}}(it)^{\mu}e^{i<x-x_{j},t>}\\
  &   & +(it)^{\mu}\cdot \overline{(it)^{\mu}}\} \hat{\psi}(t)dt\\
  & = & \frac{1}{(2\pi)^{n}}\int_{R^{n}}|\sum_{j=1}^{M}u_{j}e^{i<x_{j},t>}|^{2}\hat{\psi}(t)dt-\frac{1}{(2\pi)^{n}}\int_{R^{n}}\sum_{j=1}^{M}u_{j}\overline{(it)^{\mu}}e^{i<x_{j}-x,t>}\hat{\psi}(t)dt\\
  &   & -\frac{1}{(2\pi)^{n}}\int _{R^{n}}\sum_{j=1}^{M}\overline{u_{j}}(it)^{\mu}e^{i<x-x_{j},t>}\hat{\psi}(t)dt+(-1)^{|\mu|}\psi^{(2\mu)}(0)\\
  & = & \frac{1}{(2\pi)^{n}}\int _{R^{n}}|\sum_{j=1}^{M}u_{j}e^{i<x_{j},t>}|^{2}\hat{\psi}(t)dt-\frac{1}{(2\pi)^{n}}\int_{R^{n}}\sum_{j=1}^{M}u_{j}\overline{(it)^{\mu}}e^{i<x-x_{j},-t>}\hat{\psi}(-t)dt\\
  &   & -\frac{1}{(2\pi)^{n}}\int_{R^{n}}\sum_{j=1}^{M}\overline{u_{j}}(it)^{\mu}e^{i<x-x_{j},t>}\hat{\psi}(t)dt+(-1)^{|\mu|}\psi^{(2\mu)}(0)\\
  & = & \left\{ \begin{array}{ll}
                  U^{T}AU-2U^{T}R^{(\mu)}(x)+(-1)^{|\mu|}\psi^{(2\mu)}(0)  & \mbox{whenever $|\mu|$ is even,}\\
                  U^{T}AU+(-1)^{|\mu|}\psi^{(2\mu)}(0)                      & \mbox{whenever $|\mu|$ is odd.}
                \end{array}
        \right.
\end{eqnarray*}
Note that it's better to add a constant $(-1)^{|\mu|}$ to $\psi^{(2\mu)}(0)$ in the definition of the Kriging function.  

\bigskip

The second mistake occurs in the proof of Theorem2.3. There are three problems in the proof. The first problem is that the Fourier transform $\hat{G}_{m}(z)$ of $G_{m}$ is wrong. Let's show it.
\begin{eqnarray*}
  &   & \hat{G}_{m}(z)\\
  & = & \frac{1}{(2\pi)^{n}}\int_{R^{n}}G_{m}(t)e^{-i<z,t>}dt\\
  & = & \frac{1}{(2\pi)^{n}}\int_{R^{n}}|g_{U}(t)|^{2}e^{-\| t\| ^{2}/m^{2}}\cdot e^{-i<z,t>}dt\\
  & = & \frac{1}{(2\pi)^{n}}\int_{R^{n}}[\sum_{j=1}^{M}u_{j}e^{i<x_{j},t>}-(it)^{\mu}e^{i<x,t>}][\sum_{j=1}^{M}\overline{u_{j}}e^{-i<x_{j},t>}-\overline{(it)^{\mu}}e^{-i<x,t>}]\\
  &   & e^{-\| t\| ^{2}/m^{2}}\cdot e^{-i<z,t>}dt\\
  & = & \frac{1}{(2\pi)^{n}}\int _{R^{n}}[\sum_{j,k}u_{j}u_{k}e^{i<x_{j}-x_{k},t>}+ \sum_{j=1}^{M}u_{j}\overline{(it)^{\mu}}   e^{i<x_{j}-x,t>} \\
  &   & -\sum_{j=1}^{M}\overline{u_{j}}(it)^{\mu}e^{i<x-x_{j},t>}+(it)^{\mu}\overline{(it)^{\mu}}]e^{-\| t\| ^{2}/m^{2}}\cdot e^{-i<z,t>}dt\\
  & = & \frac{1}{(2\pi)^{n}}\int _{R^{n}}\sum_{j,k}u_{j}u_{k}e^{-\| t\| ^{2}/m^{2}}\cdot e^{-i<z-(x_{j}-x_{k}),t>}dt\\
  &   & +\frac{1}{(2\pi)^{n}}\int _{R^{n}}\sum_{j=1}^{M}u_{j}\overline{(it)^{\mu}}e^{-\| t\| ^{2}/m^{2}}\cdot e^{-i<z-(x_{j}-x),t>}dt\\
  &   & -\frac{1}{(2\pi)^{n}}\int _{R^{n}}\sum_{j=1}^{M}\overline{u_{j}}(it)^{\mu}\cdot e^{-\| t\| ^{2}/m^{2}}\cdot e^{-i<z-(x-x_{j}),t>}dt\\
  &   & +\frac{1}{(2\pi)^{n}}\int _{R^{n}}(it)^{\mu}\overline{(it)^{\mu}}e^{-\| t\| ^{2}/m^{2}}\cdot e^{-i<z,t>}dt\\
  & = & \sigma_{mn}m^{n}\sum_{j,k}u_{j}u_{k}e^{-\| z-(x_{j}-x_{k})\| ^{2}m^{2}/4}\\
  &   & +\sigma_{mn}m^{n}\sum_{j=1}^{M}u_{j}D^{\mu}e^{-\| z-(x_{j}-x)\| ^{2}m^{2}/4}\\
  &   & \sigma_{mn}m^{n}(-1)^{|\mu|}\sum_{j=1}^{M}u_{j}D^{\mu}e^{-\| z-(x-x_{j})\| ^{2}m^{2}/4}\\
  &   & +\sigma_{mn}m^{n}(-1)^{|\mu|}D^{2\mu}e^{-\| z\| ^{2}m^{2}/4}
\end{eqnarray*}
, where $\sigma_{mn}$ is some constant independent of $m$.

\bigskip

Note that the second and third items of the result cannot be incorporated whenever $|\mu|$ is odd. Now, the second problem arises when we insert $\hat{G}_{m}$ into (17). We display the process.

\begin{eqnarray}
  &   & lim_{m\rightarrow \infty}\frac{1}{(2\pi)^{n}}\int_{R^{n}}\hat{G}_{m}(z)\psi(z)dz \nonumber \\
  & = & lim_{m\rightarrow \infty} \frac{1}{(2\pi)^{n}}\cdot \sigma_{mn}\sum_{j,k}u_{j}u_{k}\int_{R^{n}}m^{n}e^{-\| z-(x_{j}-x_{k})\| ^{2}m^{2}/4}\cdot \psi(z)dz \nonumber \\
  &   & +lim_{m\rightarrow \infty }\frac{1}{(2\pi)^{n}}\cdot \sigma_{mn}\sum_{j=1}^{M}\int_{R^{n}}m^{n}D^{\mu}e^{-\| z-(x_{j}-x)\| ^{2}m^{2}/4}\cdot \psi(z)dz \nonumber \\
  &   & +lim_{m\rightarrow \infty }\frac{1}{(2\pi)^{n}}\cdot \sigma_{mn}(-1)^{|\mu|}\sum_{j=1}^{M}\int_{R^{n}}m^{n}D^{\mu}e^{-\| z-(x-x_{j})\| ^{2}m^{2}/4}\cdot \psi(z)dz \nonumber \\
  &   & +lim_{m\rightarrow \infty}\frac{1}{(2\pi)^{n}}\cdot \sigma_{mn}(-1)^{|\mu|}\int_{R^{n}}m^{n}D^{2\mu}e^{-\| z\| ^{2}m^{2}/4}\cdot \psi(z)dz
\end{eqnarray}
Wu and Schaback conclude that the limits of the four integrals of the right-hand side are equal to
$$\psi(x_{j}-x_{k}),\ (-1)^{|\mu|}D^{\mu}\psi(x_{j}-x),\ (-1)^{|\mu|}D^{\mu}\psi(x-x_{j})\ and\ (-1)^{|2\mu|}D^{2\mu}\psi(0)$$
respectively, by applying $\delta$-function. However, by p.217 of \cite{J}, $\psi$ should be a good function when applying $\delta$-function. That is, $\psi$ should be a member of the Schwartz space. Fortunately, as long as we require $\psi$ to be in $C^{1}$, then $\delta$-function still works even if $\psi$ is not a good function. This can be shown as follows.
$$\int_{-\infty}^{\infty}(m/\pi)^{n/2}e^{-mr^{2}}\psi(x)dx=\psi(0)+\int_{-\infty}^{\infty}(m/\pi)^{n/2}\{ \psi(x)-\psi(0)\} e^{-mr^{2}}dx.$$
Since
\begin{eqnarray*}
  |\psi(x)-\psi(0)| & = & |\int_{0}^{x_{1}}\frac{\partial \psi(t_{1},0,\ldots ,0)}{\partial t_{1}}dt_{1}+\int_{0}^{x_{2}}\frac{\partial \psi(x_{1},t_{2},0,\ldots ,0)}{\partial t_{2}}dt_{2}+ \cdots \\
                    &   & + \int_{0}^{x_{n}}\frac{\partial \psi(x_{1},x_{2},\ldots ,t_{n})}{\partial t_{n}}dt_{n}|\\
                    & \leq & M(|x_{1}|+|x_{2}|+\cdots +|x_{n}|)
\end{eqnarray*}
, the integral does not exceed $Mn/m^{1/2}$ which tends to zero as $m\rightarrow \infty $ and 
$$\int_{-\infty }^{\infty}\delta(x)\psi(x)dx=\psi(0)$$
is justified.

\bigskip

What's making trouble is the derivative. This is the third problem. We know that for any generalized function $g$ and good function $\gamma$, 
$$\int_{R^{n}}\{ \partial_{1}^{q_{1}}\cdots \partial_{n}^{q_{n}}g(x)\} \gamma(x)dx=(-1)^{q}\int_{R^{n}}g(x)\partial^{q}\gamma(x)dx$$
where $q=(q_{1},\ldots ,q_{n})$. The requirement that $\gamma(x)$ be a good function is absolutely necessary and cannot be modified. Thus the second and third integrals in the right-hand side of (25) cannot be reduced to $(-1)^{|\mu|}D^{\mu}\psi(x_{j}-x)$ and $(-1)^{|\mu|}D^{\mu}\psi(x-x_{j})$ respectively. Moreover, the item $\psi^{(2\mu)}(0)$ needs a coefficient $(-1)^{|\mu|}$ also. These bad mistakes destroy (\ref{2K}) thoroughly.

\bigskip

The third crucial mistake occurs in (22) whenever $\hat{f}$ and $\hat{\psi}$ are generalized functions. Note that (22) is a result of (\ref{2M}) which is also proven as follows for the generalized case.

\bigskip

Let $g_{U}(t)=\sum_{j=1}^{M}u_{j}^{(\mu)}(x)e^{i<x_{j},t>}-(it)^{\mu}e^{i<x,t>}$. Then, by Lebesgue Convergence Theorem, 
\begin{eqnarray*}
  \frac{1}{(2\pi)^{n}}\int_{R^{n}}g_{U}(t)\hat{f}(t)dt & = & lim_{m\rightarrow \infty}\frac{1}{(2\pi)^{n}}\int_{R^{n}}g_{U}(t)e^{-\| t\| ^{2}/m^{2}}\hat{f}(t)dt\\
                                                       & = & lim_{m\rightarrow \infty}\frac{1}{(2\pi)^{n}}\int_{R^{n}}\hat{G}_{m}(z)f(z)dz\\
                                                       &   & where\ G_{m}(t):=g_{U}(t)e^{-\| t\| ^{2}/m^{2}}.
\end{eqnarray*}
The Fourier transform $\hat{G}_{m}(z)$ can be calculated as follows.

\begin{eqnarray*}
  \hat{G}_{m}(z) & = & \frac{1}{(2\pi)^{n}}\int_{R^{n}}[\sum_{j=1}^{M}u_{j}^{(\mu)}(x)e^{i<x_{j},t>}-(it)^{\mu}e^{i<x,t>}]e^{-\| t\| ^{2}/m^{2}}\dot e^{-i<z,t>}dt\\
                 & = & \frac{1}{(2\pi)^{n}}\sum_{j=1}^{M}u_{j}^{(\mu)}(x)\int_{R^{n}}e^{-\| t\| ^{2}/m^{2}}\cdot e^{-i<z-x_{j},t>}dt\\
                 &   & -\frac{1}{(2\pi)^{n}}\int_{R^{n}}(it)^{\mu}e^{-\| t\| ^{2}/m^{2}}\cdot e^{-i<z-x,t>}dt\\
                 & = & \frac{1}{(2\pi)^{n}}m^{n}\sum_{j=1}^{M}u_{j}^{(\mu)}(x)\cdot e^{-\frac{m^{2}}{4}\| z-x_{j}\| ^{2}} -\frac{1}{(2\pi)^{n}}m^{n}(-1)^{|\mu|}D^{\mu}\int_{R^{n}}e^{-\frac{\| t\| ^{2}}{m^{2}}}\cdot e^{-i<z-x,t>}dt\\
                 & = & \frac{1}{(2\pi)^{n}}m^{n}\sum_{j=1}^{M}u_{j}^{(\mu)}(x)\cdot e^{-\frac{m^{2}}{4}\| z-x_{j}\| ^{2}}-\frac{1}{(2\pi)^{n}}m^{n}(-1)^{|\mu|}D_{z}^{\mu}e^{-\frac{m^{2}}{4}\| z-x\| }.
\end{eqnarray*}
Now, 
\begin{eqnarray*}
  &   & lim_{m\rightarrow \infty}\frac{1}{(2\pi)^{n}}\int_{R^{n}}\hat{G}_{m}(z)f(z)dz\\
  & = & lim_{m\rightarrow \infty}\sigma_{mn}[1/(2\pi)^{n}]^{2}\int_{R^{n}}m^{n}\sum_{j=1}^{M}u_{j}^{(\mu)}(x)e^{-\frac{m^{2}}{4}\| z-x_{j}\| ^{2}}\cdot f(z)-(-1)^{|\mu|}m^{n}D^{\mu}_{z}e^{-\frac{m^{2}}{4}\| z-x\| ^{2}}\cdot f(z)dz\\
  & = & \frac{\sigma_{mn}}{(2\pi)^{2n}}\sum_{j=1}^{M}u_{j}^{(\mu)}(x)lim_{m\rightarrow \infty}\int_{R^{n}}m^{n}e^{-\frac{m^{2}}{4}\| z-x_{j}\| ^{2}}\cdot f(z)dz\\
  &   & -\frac{\sigma_{mn}}{(2\pi)^{2n}}\cdot (-1)^{|\mu|}lim_{m\rightarrow \infty}\int_{R^{n}}m^{n}D_{z}^{\mu}e^{-\frac{m^{2}}{4}\| z-x\| ^{2}}\cdot f(z)dz\\
  & = & \frac{\sigma_{mn}}{(2\pi)^{3n/2}}\sum_{j=1}^{M}u_{j}^{(\mu)}(x)f(x_{j})-\frac{\sigma_{mn}}{(2\pi)^{3n/2}}(-1)^{|\mu|}\cdot (-1)^{|\mu|}\cdot  D^{\mu}f(x)\\
  & = & \frac{\sigma_{mn}}{(2\pi)^{3n/2}}[s^{(\mu)}(x)-f^{(\mu)}(x)].
\end{eqnarray*}
Hence we get (22) up to a constant for the generalized case. 

\bigskip

In the above course, we also use the technique of $\delta$-function and the derivative of generalized function. Mistakes occur again since $f$ is not a good function. However, if $|\mu|=0$, the above argument still holds as long as we require $f$ to be a $C^{1}$-function.

\bigskip

We conclude that Wu and Schaback's error bound can be used only when $|\mu|=0$ with an additional restriction that both $\psi$ and $f$ should be in $C^{1}$ classically. What's noteworthy is that one should avoid using this error bound when $\psi$ is a thin-plate spline. As for multiquadrics and Gaussians, it still works.

\bigskip

The significance of derivatives should be explained here. We know that derivatives reflect smoothness. If the approximand and approximant can match in smoothness, it will be nice. I fear  r.b.f. people probably overemphasize the infinitely differentiable multiquadrics. Of course, multiquadrics are a very useful tool in the multivariate scattered data interpolation. However, as pointed out by Yoon in \cite{Yo1}, most data points do not arise from extremely smooth functions. If this is true, multiquadrics are not very satisfactory in most cases. An example can be given here. It's well known that Weierstrass raised a function which is everywhere continuous, but is everywhere undifferentiable. If we approximate such a function by multiquadrics, it's meaningless even though its error bound is small. Although this example is a bit too drastic, it provides an evidence that sometimes we do have to consider the matching of smoothness. What's noteworthy is that in the practical problems of differential equations, most solutions are at most $C^{2}$ or $C^{3}$ functions. They are not very smooth.

\bigskip

I personally believe that unsmooth functions should be approximated by unsmooth functions. In fact, this is the central idea of this paper.

\section{Comparison of the function spaces}

In Wu and Schaback's approach, the space of approximated functions $F_{\phi}$ is defined to be the set of all functions from $R^{n}$ into $R$ which has a generalized Fourier transform $\hat{f}$ coinciding on $R^{n}-\{ 0\} $ with a continuous function satisfying $$\int_{R^{n}}|\hat{f}(t)|^{2}(\hat{\psi}(t))^{-1}dt<\infty $$
for $\psi(x)=\phi(\| x\| )$ satisfying (\ref{2I}) and (\ref{2J}).

\bigskip

Many people mistake Wu and Schaback's space for Madych and Nelson's space of approximands. This may lead to a disaster.

\bigskip

Madych and Nelson's space \cite{MN1},\cite{MN2} is much more complicated to describe than that of Wu and Schaback. However, L-T. Luh has made a clear characterization for Madych and Nelson's space which is easier to understand. Let's follow Luh's notation and definitions in \cite{Lu1} and \cite{Lu2}.

\bigskip

In the following definitions, $\Omega$ denotes an arbitrary subset of $R^{n}$ and $P_{m}^{n}$ denotes the set of $n$-variable polynomials of order not exceeding $m$.
\begin{definition}
  $(P_{m}^{n})_{\Omega}^{\perp}:=\{ \sum_{i=1}^{N}c_{i}\delta_{x_{i}}| c_{i}\in C,\ x_{i}\in \Omega\ for\ 1\leq i\leq N,\ \sum_{i=1}^{N}c_{i}p(x_{i})
 \  \ \ \ \ \ \ \ \ \ \ \\ \ \ \ \ \ \  =0\ for\ all\ p\in P_{m}^{n}\} $ where $\delta_{x_{i}}$ is a measure defined by
$\delta_{x_{i}}(E)=\left\{   \begin{array}{lll}
                            1 & if & x_{i}\in E \\
                            0 & if & x_{i}\not\in E
                          \end{array} \right. $
                    for any subset $E$ of $R^{n}$.
\end{definition} 
We can construct a norm in $(P_{m}^{n})_{\Omega}^{\perp}$. If we define 
$$(\lambda,\mu)_{\psi}:=\sum_{i=1}^{N}\sum_{j=1}^{M}\lambda_{i}\overline{\mu}_{j}\psi(x_{i}-y_{j})$$
for $\lambda,\mu\in (P_{m}^{n})_{\Omega}^{\perp}$ with $\lambda=\sum_{i=1}^{N}\lambda_{i}\delta_{x_{i}},\ \mu=\sum_{j=1}^{M}\mu_{j}\delta_{y_{j}}$, it's easily seen $(\cdot ,\cdot )_{\psi}$ is an inner product and induces a norm in $(P_{m}^{n})_{\Omega}^{\perp}$ as long as $\psi$ is $c.p.d.$ of order $m$ in $\Omega$.  
\begin{definition}
  $\overline{(P_{m}^{n})_{\Omega}^{\perp}}$ is the Hilbert space completion of $(P_{m}^{n})_{\Omega}^{\perp}$ under the norm induced by $\psi$.
\end{definition}

\begin{definition}
  $\underline{F_{\phi,\Omega}}:= \{ f(x)=(\mu,\delta_{(x)})_{\psi}| \mu \in \overline{(P_{m}^{n})_{\Omega}^{\perp}},\ x \in    \Omega,\ (\cdot ,\cdot )_{\psi}\ is\ the\ inner  $\ \ \ \ \ \ \ \ \ \ \ \ \ \ \ \ \ \ \ \ \ \ \ \ \ \ \ \ \ \ \ \ \ \ \ \ \ \ \ \ \ \ \ \ \ \ \ \ \ \ \ \ \ \ \ \ \ \ \ \ \ \ \ \ \ \ $  product\ in\ \overline{(P_{m}^{n})_{\Omega}^{\perp}}\ induced\ by\ \psi \} $ where $\delta_{(x)}:=\delta_{x}-\sum_{i=1}^{m'}l_{i}(x)\delta_{x_{i}}$ with $l_{1},\ldots ,l_{m'}$ a basis of $P_{m}^{n}$ and $x_{1},\ldots ,x_{m'}\in \Omega $ satisfying $l_{i}(x_{j})=\delta_{ij}$.
\end{definition}

\begin{definition}
  Let $\psi$ be $c.p.d.$ of order $m\geq 0$ in $\Omega$. Then 
$$C_{\psi}(\Omega)=P_{m}^{n}(\Omega)\oplus \underline{F_{\psi,\Omega}}$$
where $\oplus $ denotes the direct sum of two linear spaces.
\end{definition}

\bigskip

In fact, Madych and Nelson's space of approximands is just $C_{\psi}(R^{n})$.

\bigskip

In \cite{WS} Wu and Schaback say that the completion of $F_{\phi}$ is just $\underline{F_{\psi,R^{d}}}$. Some authors even say that $F_{\phi}$ is equal to $C_{\psi}(R^{d})$. See \cite{Yo1} and \cite{Yo2}.
\bigskip

In order to make a lucid clarification, I must quote two theorems of \cite{MN1}.
\begin{theorem}
  Let $\psi$ be continuous and conditionally positive definite of order $m$. Then it is possible to choose a positive Borel measure $\mu$ on $R^{n}-\{ 0\} $, constants $a_{\gamma},\ |\gamma |\leq 2m$ and a function $\chi $ in $\cal D$ such that: $1-\hat{\chi}(\xi)$ has a zero of order $2m+1$ at $\xi=0$; both of the integrals $\int_{0<|\xi|<1}|\xi|^{2m}d\mu(\xi)$, $\int_{|\xi|\geq 1}d\mu(\xi)$ are finite; for all $\sigma\in \cal D$,
\end{theorem}
\begin{equation}
  \int \psi(x)\sigma(x)dx = \int [\hat{\sigma}(\xi)-\hat{\chi}(\xi)\sum_{|\gamma|<2m}D^{\gamma}\hat{\sigma}(0)\frac{\xi ^{\gamma}}{\gamma !}]d\mu (\xi)+ \sum _{|\gamma|\leq 2m}D^{\gamma}\hat{\sigma}(0)\frac{a_{\gamma}}{\gamma !}. \label{4A}
\end{equation}
This uniquely determines the measure $\mu$ and the constants $a_{\gamma}$ for $|\gamma |=2m$. In addition, for every choice of complex numbers $c_{\alpha},|\alpha|=m$, 
$$\sum_{|\alpha|=m}\sum_{|\beta|=m}a_{\alpha +\beta}c_{\alpha}\overline{c}_{\beta}\geq 0.$$
\begin{theorem}
  Let $m,\psi,\mu$ and $a_{\gamma}$ be as in Theorem4.1. Assume $d\mu(\xi)=w(\xi)d\xi $ and $a_{\gamma}=0$ for all $|\gamma|=2m$. Let $\rho$ be the Borel measure on $R^{n}$ defined by $d\rho(\xi)=r(\xi)d\xi$, where $r$ is defined by 
$$r(\xi)=\frac{1}{(2\pi)^{2n}|\xi|^{2m}w(-\xi)},$$
with $r(\xi)=\infty$ when $w(-\xi)=0$. Then $f\in C_{\psi}(R^{n})$ if and only if $f\in \cal S'$ and $(D^{\alpha}f)^{\wedge }\in L^{2}(\rho)$ for every $|\alpha|=m$. In that case, $(f,f)_{\psi}$ is given by 
$$\sum_{|\alpha|=m}\frac{m!}{\alpha !}\| (D^{\alpha}f)^{\wedge }\|^{2}_{L^{2}(\rho)}=\int |g|^{2}d\mu=(f,f)_{\psi}$$
where $g$ is defined in Proposition3.2 of \cite{MN2}. 
\end{theorem}
By Theorem4.2, one easily finds that Wu and Schaback's $F_{\phi}$ is equal to a subset of  Madych and Nelson's $C_{\psi}(R^{n})$  when $m=0$, with the same norm. Here, one should be careful. Wu and Schaback's space consists only of real-valued functions. However, Madych and Nelson's space contains complex-valued functions. When $m>0$, the two spaces cannot be compared. For $m>0$, the Hilbert space completion of $F_{\phi}$ still has nothing to do with $C_{\psi}(R^{n})$.

\bigskip

Now, another big problem of Wu and Schaback's error bound can be pointed out. That is, there is no reason to believe that the generalized derivative $f^{(\mu)}(x)$ in (\ref{2P}) coincides with a nicely behaving classical function. This question is absolutely nontrivial. Although some people may add this condition as an additional requirement, this will become a drawback because it's rather restrictive. Madych and Nelson deal with this problem perfectly in their paper \cite{MN2}. In their theory the generalized derivative coincides naturally with a nicely behaving  classical function.However, it's suitable only for $C_{\psi}(R^{n})$. Moreover, Madych and Nelson's space $C_{\psi}(R^{n})$ seems to be more meaningful than Wu and Schaback's space $F_{\phi}$ for two reasons. First, it has a natural interpretation in physics. Second, it contains a dense subset which consists of natural extensions of the univariate cubic spline. The two reasons are very strong. All these can be found in \cite{MN1}.

\bigskip

In fact, not only Wu, Schaback and Yoon, but also many others misunderstand the two spaces. The main reason is that Madych and Nelson's papers are very difficult to read. In order to get a better understanding to their papers, it might be helpful to read \cite{Lu1} and \cite{Lu2} first.

\bigskip

Finally, I would like to point out a key point of the comparison of the two spaces and the two error bounds. It's well known that there are two notions of generalized Fourier transform. One is the Russian notion of Gelfand and Schilov. The other is the French notion of Schwartz. Many people think the notion used by Madych and Nelson in \cite{MN2} is the Russian notion. In fact, this is wrong. What they use is still the French notion. This can be seen in page 213 and page 219 of \cite{MN2}. Madych and Nelson studied Duchon's works and tried to promote it. That probably is the reason why they use the French notion of generalized Fourier transform. Unfortunately, in the bibliography of \cite{MN2}, only Gelfand and Vilenkin's book "Generalized Functions", Vol.4, 1964 appears. None of Schwartz, Jones, or Stein/Weiss's book appears. This misleads people into thinking that the notion used by them is the Gelfand/Schilov's definition. This seems to be part of the reasons many people do not like Madych and Nelson's approach. They prefer Wu and Schaback's approach which uses the better digestible French notion.

 Department of Applied Mathematics,

Providence University,

Shalu Town,

Taichung County,

Taiwan

\bigskip

Email:ltluh@pu.edu.tw

\end{document}